\documentclass[12pt]{amsart}

\newtheorem{theorem}{Theorem}[section]
\newtheorem{proposition}{Proposition}[section]
\newtheorem{lemma}[theorem]{Lemma}

\def\C{{\mbox{\rm\kern.24em
\vrule width.03em height1.43ex depth-.052ex \kern-.26em C}}}
\def\QSet{\mbox{\rm\kern.24em
\vrule width.03em height1.48ex depth-.051ex \kern-.26em Q}}
\def\Z{{\mbox{\rm\kern.25em
\vrule width.03em height0.57ex depth0ex
\kern.033em
\vrule width.03em height1.52ex depth-0.96ex \kern-.338em Z}}}

\def\N{{{\mbox{\rm I\kern-.2em N}}_0}}
\def\M{{{\mbox{\rm I\kern-.2em M}}}}
\def\PSet{\mbox{\rm I\kern-.22em P}}
\def\R{{\mbox{\rm I\kern-.22em R}}}

\def\P{{\bf P}}
\def\p{{\bf p}}

\def\I{{\bf I}}

\def\T{{\bf T}}

\def\I{{\bf I}}

\def\<{\left<}
\def\>{\right>}

\def\size{{\bf size}}
\def\w{{\tilde{w}}}

\title[New uniform bounds for a Walsh model]{New uniform bounds for a Walsh model of the bilinear Hilbert transform}

\author
[R. Oberlin, C. Thiele]
{Richard Oberlin \ \ \ Christoph Thiele }

\address{R. Oberlin, Department of Mathematics,
UCLA, Los Angeles, CA 90095-1555, USA}
\email{oberlin@math.ucla.edu}

\address{C. Thiele, Department of Mathematics,
UCLA, Los Angeles, CA 90095-1555, USA}
\email{thiele@math.ucla.edu}


\thanks{R.O. partially supported by NSF VIGRE grant
DMS 0502315.
C.Th. partially supported by NSF grant DMS 0701302.}

\date{\today}

\begin{document}

\begin{abstract}
We prove old and new $L^p$ bounds for the quartile operator, a Walsh 
model of the bilinear Hilbert transform, uniformly in the parameter that 
models degeneration of the bilinear Hilbert transform. We obtain the full
range of exponents that can be expected from known bounds in 
the degenerate and non-degenerate cases. For the new estimates with 
exponents $p$ 
close to $1$ the argument relies on a multi-frequency Calderon-Zygmund 
decomposition.
\end{abstract}

\maketitle

\section{Introduction}

The notion of a bilinear Hilbert transform usually refers to a member of
a family of bilinear operators parameterized by a unit vector $\beta$ perpendicular 
to $(1,1,1)$. We will write the bilinear operators in this family more symmetrically as
dual trilinear forms $\Lambda_\beta$, acting on three test functions on the real line:
$$\Lambda_\beta(f_1,f_2,f_3):=\int p.v. \int 
f_1(x-\beta_1 t)f_2(x-\beta_2 t)f_3(x-\beta_3 t)\, \frac{dt}t\ dx\ \ .$$
The interesting case, which we call non-degenerate, is when the
three components of $\beta$ are pairwise different.
If two of the components of $\beta$ are equal the form reduces to the combination 
of a pointwise product and the dual
of the classical linear Hilbert transform.
A priori $L^p$ bounds in the non-degenerate case were first shown in 
\cite{laceyt1} and \cite{laceyt2}. Namely,
for each $1<p_1,p_2,p_3\le \infty$ with $\sum_j 1/p_j=1$ we have
\begin{equation}\label{pbounds}
|\Lambda_\beta(f_1,f_2,f_3)|\le C_{\beta,p_1,p_2,p_3} \prod_{j=1}^3 \|f_j\|_{p_j}\ \ .
\end{equation}
The condition $\sum_j 1/p_j=1$ is necessary by dilation symmetry of the
form $\Lambda_\beta$ and shall be assumed throughout the rest of this discussion.

If each $f_j$ is bounded by the characteristic function of a set $E_j$, then 
inequality (\ref{pbounds}) implies the restricted type estimate 
\begin{equation}\label{alphabounds}
|\Lambda_\beta(f_1,f_2,f_3)|\le C_{\beta,\alpha_1,\alpha_2,\alpha_3} \prod_{j=1}^3 |E_j|^{\alpha_j}
\end{equation}
where $\alpha_j=1/p_j$ satisfies $0\le \alpha_j<1$.
More generally, the argument in \cite{laceyt2} shows that inequality (\ref{alphabounds})
continues to hold in the range $-1/2<\alpha_j <1$ under the additional
assumption that if $\alpha_j<0$ then $f_j$ is bounded by the characteristic function of a 
major subset ${E_j}'\subset E_j$ that depends on the sets $E_1$, $E_2$, and $E_3$. Here a major
subset is one of measure at least half the measure of the ambient set.
The passage to a major subset of $E_j$ is natural and necessary in the setting of negative 
exponents and was introduced in this context in \cite{mtt0}.

The range of triples $(\alpha_1,\alpha_2,\alpha_3)$ for which one has
the a priori estimate (\ref{alphabounds}) appears in Figure \ref{hexagon}
as the the convex hull of the open triangles $a_1$, $a_2$, and $a_3$.
Note that the closed triangle $c$ represents the
local $L^2$ case with $2\le p_1,p_2,p_3\le \infty$, while
the convex  hull of the open triangles $b_1,b_2,b_3$ represents
the reflexive Banach triangle where $1<p_1,p_2,p_3<\infty$.

\noindent
\setlength{\unitlength}{0.9mm}
\begin{figure}\label{hexagon}
\begin{picture}(72,78)
\put(46,60){\circle*{1.6}}
\put(28,30){\circle*{1.6}}
\put(64,30){\circle*{1.6}}
\put(46,30){\circle*{0.8}}
\put(19,45){\circle*{0.8}}
\put(28,60){\circle*{0.8}}
\put(64,60){\circle*{0.8}}
\put(73,45){\circle*{0.8}}
\put(37,15){\circle*{0.8}}
\put(55,15){\circle*{0.8}}
\put(37,45){\circle*{0.8}}
\put(55,45){\circle*{0.8}}

\put(10,28){$(0,0,1)$}
\put(68,28){$(0,1,0)$}
\put(40,64){$(1,0,0)$}

\put(27,49){$a_2$}
\put(63,49){$a_3$}
\put(45,38){$c$}
\put(45,19){$a_1$}
\put(45,49){$b_1$}
\put(54,34){$b_2$}
\put(36,34){$b_3$}

\put(25,38){$d_{23}$}
\put(61,38){$d_{32}$}
\put(34,53){$d_{21}$}
\put(52,53){$d_{31}$}
\put(34,23){$d_{13}$}
\put(52,23){$d_{12}$}

\put(10,60){\line(1,0){72}}
\put(10,60){\line(3,-5){36}}
\put(82,60){\line(-3,-5){36}}
\put(19,45){\line(3,5){9}}
\put(73,45){\line(-3,5){9}}
\put(37,15){\line(1,0){18}}

\put(19.4,45){\line(1,0){53.2}}
\put(28.4,30){\line(1,0){35.2}}
\put(28.2,59.67){\line(3,-5){26.6}}
\put(46.2,59.67){\line(3,-5){17.6}}
\put(28.2,30.33){\line(3,5){17.6}}
\put(37.2,15.33){\line(3,5){26.6}}





\end{picture}
\end{figure}

\noindent

In the degenerate case, say $\beta_2=\beta_3$, a priori estimates
estimates for $\Lambda_\beta$ 
follow from H\"older's inequality and bounds
for the linear Hilbert transform. One has bounds of the  type (\ref{alphabounds}) 
if $\alpha_1>0$ and $\alpha_1,\alpha_2,\alpha_3<1$. The intersection of this region with the 
region of bounds for the non-degenerate case is the convex hull
of the open triangles $a_2,a_3,b_2,b_3$. It is natural to
ask whether one has bounds for the non-degenerate case uniformly
in the parameter $\beta$ in a small neighborhood of the degenerate 
case $\beta_2=\beta_3$. Several articles have been written on this question: 
\cite{thiele4} proves inequality (\ref{alphabounds}) uniformly in such $\beta$ 
at the two upper corners of the triangle $c$ under the assumption $f_j$ is supported 
on a major subset ${E_j}'\subset E_j$ when $\alpha_j=0$. Grafakos and Li \cite{grafakosl1}
show inequality (\ref{pbounds}) in the triangle $c$ and Li \cite{li} shows (\ref{alphabounds})
uniformly in the open triangles $a_2$ and $a_3$. One can interpolate these results
to get bounds in the convex hull of the open triangles $a_2$, $a_3$, and $c$, but 
it remains open to date whether uniform bounds hold in the entire open triangles $b_2$ and $b_3$. 
The current paper presents progress in this direction by proving uniform bounds in
$b_2$ and $b_3$ for a discrete model of the bilinear Hilbert transform. 

The quartile operator was introduced in \cite{thiele} as a discrete model
for the non-degenerate bilinear Hilbert transform. In \cite{thiele5} a family of 
related operators was introduced that models the set of Hilbert transforms near the 
degenerate case and allows to address uniformity questions in the model case. 
Moreover, inequality (\ref{alphabounds}) was shown at the two
upper corners of the triangle c under the assumption $f_j$ is supported on a major
subset of $E_j$ if $\alpha_j=0$. In the current paper, we extend these results to the 
entire convex hull of the open triangles $a_2,a_3,b_2,b_3$ and thus the full range in which 
we know bounds both for the degenerate and the non-degenerate case. 
Our proof simplifies that in \cite{thiele4}, using an  approach via phase plane 
projections developed in the continuous case in \cite{mtt4}. It also 
uses a simple discrete version of the multi-frequency Calderon Zygmund 
decomposition introduced in \cite{not} as well as a technique of \cite{osttw} of using 
BMO bounds for the counting function defined further below. In this sense this article 
also serves as expository survey of these techniques in the discrete setting. We plan to
address the extension of the novel results in this paper to the continuous setting
in future work.

We proceed to formulate the main theorems of this paper in detail. 
The Walsh phase plane is the closed first quadrant $\R_+\times \R_+$ 
of the plane.  A dyadic rectangle is a rectangle in the Walsh phase plane 
of the form 
\begin{equation}\label{tiledef}
p=I\times \omega= [2^kn,2^k(n+1))\times [2^{-k'} l, 2^{-k'} (l+1))
\end{equation}  with integers $k,k',n,l$ 
and $0\le n,l$. A tile is a dyadic rectangle of area one, while a bitile
is a dyadic rectangle of area two. Each bitile can be split into
upper tile $P_u$ and lower tile $P_d$, or alternatively into
left tile $P_{\rm left}$ and right tile $P_{\rm right}$. Associated to each tile $p$ 
is a Walsh wave packet $w_p$, which is a certain function in $L^2(\R_+)$
normalized to have $L^2$ norm one. We will also use the abbreviation $\w_p = |I_p|^{1/2} w_p$ for the $L^\infty$ normalized wave-packet.
With the notation as in (\ref{tiledef}), if $l=0$, then this
wave packet is defined as the appropriate multiple of the characteristic 
function of $I$. For other values of $l$ 
it is defined recursively via the identities
$$w_{P_u}= ( w_{P_{\rm left}}-w_{P_{\rm right}})/\sqrt{2}\ \ ,$$
$$w_{P_d}= ( w_{P_{\rm left}}+w_{P_{\rm right}})/\sqrt{2}\ \ .$$
By induction on the depth of this recursion one can show 
(\cite{thiele}) that $w_p$ is supported on $I$, it has
constant modulus on $I$, and disjoint tiles correspond to 
orthogonal wave packets. If $S$ is a subset of the Walsh 
phase plane that can be written as a disjoint union of a collection
$\p$ of tiles, we define the phase plane projection
associated to $S$ to be the orthogonal projection 
$$\Pi_S f =\sum_{p\in \p} \<f,w_p\>w_p\ \ .$$
One can show that this projection is independent 
of the particular tiling $\p$ of the set $S$, justifying the
notation that ignores the particular choice of tiling.
For a subset $S$ in the phase plane and an integer $L$ 
define $2^LS$ to be the set $\{(x,2^L\xi), (x,\xi)\in S\}$.

We define the quartile\footnote {A quartile is a dyadic rectangle of area four. The name
quartile form is inherited from the use of quartiles in \cite{thiele} to define
a related model.} 
form with parameter $L\ge 2$ as follows:
$$\Lambda_L(f_1,f_2,f_3)
:=\int \sum_{P} \w_{P_d}(x) \Pi_{P_u} f_1(x)
\prod_{j=2}^3 \Pi_{2^LP_d} f_j(x) \, dx\ \ .$$
Here $P$ runs through the set of all bitiles. To avoid technical arguments
we shall restrict this set to the set of all bitiles contained in the strip
$\R_+\times [0,2^N)$ for some very large $N$. This restriction is equivalent 
to assuming that $f_1$ is constant on intervals of length $2^{-N}$.
The bounds claimed in the following theorems are independent of $N$. While fixing 
$N$ destroys the dilation symmetry of the form $\Lambda_L$, the family of $\Lambda_L$ 
for all such $N$ retains the dilation symmetry and so do the main theorems.
We will avoid explicit mentioning of $N$ in most of this paper.

Our main results are
the following two theorems:
\begin{theorem}\label{main1}
For any exponents $1<p_1,p_2,p_3<\infty$ with
$$1/p_1+1/p_2+1/p_3=1$$
there is a constant $C_{p_1,p_2,p_3}$ independent of $L$ (and $N$) such that we have the
a priori estimate
$$|\Lambda_L(f_1,f_2,f_3)|\le C_{p_1,p_2,p_3}\prod_{j=1}^3 \|f_j\|_{p_j}\ \ .$$
\end{theorem}
\begin{theorem}\label{main2}
Let $0<\alpha_1,\alpha_3<1$ and $-1/2<\alpha_2\le 0$
with $\sum_j\alpha_j=1$. For any three measurable subsets  $E_j$, $j=1,2,3$ 
of $\R_+$ such that $|E_2|$ is maximal among the $|E_j|$ there is a major 
subset ${E_2}'$ of ${E_2}$ such that for any three measurable functions $f_j$,
bounded in absolute value by the characteristic function of $E_j$ if $j\neq2$ and
the characteristic function of ${E_2}'$ if $j=2$, we have the following estimate
$$|\Lambda_L(f_1,f_2,f_3)|\le C_{\alpha_1,\alpha_2,\alpha_3}
|E_1|^{\alpha_1} |E_2|^{\alpha_2}|E_3|^{\alpha_3}$$
uniformly in the parameter $L$ (and $N$).
\end{theorem}

Note that if $|E_2|$ is not maximal among the $|E_j|$, the conclusion of
Theorem \ref{main2} follows with ${E_2}'=|E_2|$ from an application of Theorem \ref{main1}
with a different set of exponents. Since the quartile form is symmetric in the 
indices $j=2,3$, one obtains 
as corollary a symmetric version of Theorem \ref{main2}. The proofs of
these theorems are sufficiently robust to allow for a perturbation of the
quartile form by an arbitrary bounded sequence
$|c_P|\le 1$:
$$\Lambda_L(f_1,f_2,f_3)
:=\int \sum_{P} \left[c_P \w_{P_d}(x)\right] \Pi_{P_u} f_1(x)
\prod_{j=2}^3 \Pi_{2^LP_d} f_j(x) \, dx\ \ .$$
This flexibility adds to the usefulness of our arguments as a model situation
for bilinear singular integrals.

In Section \ref{trees} we prove estimates for trees,  which are collections 
of bitiles with lacunary structure, and present a tree selection algorithm.
In Section \ref{proof} we use these ingredients to assemble the proofs of
Theorems \ref{main1} and \ref{main2}.

\section{Trees} \label{trees}

A dyadic rectangle $P=I\times \omega$ is less than or equal to another dyadic
rectangle $P'=I'\times \omega'$,
in writing $P\le P'$, if $I\subset I'$ and $\omega'\subset \omega$. 
Two dyadic rectangles of the same area are comparable under this
order relation if and only if they have nonempty intersection. 

A tree $T$ is a collection of bitiles with a unique maximal element, 
usually denoted by $P_T=I_T\times \omega_T$.

A set $\P$ of bitiles is called convex if for any two 
elements $P_1,P_2\in \P$ and any bitile $P$ which satisfies $P_1<P<P_2$ we have 
$P\in \P$. It is shown in \cite{thiele} by induction on the number of bitiles
that for any convex set $\P$ of bitiles the union $\bigcup_{P\in \P} P$ may 
be written as the disjoint union of tiles. If the set $\P$ is a convex 
tree $T$, one such  tiling is obtained by
decomposing the tree as union of three collections of bitiles 
\begin{equation}\label{upperlowertree}
 T=\{P_T\}\cup T_u\cup T_d
\end{equation}
where
$$T_d=\{P\in T: P_T\le P_d\}\ \ ,$$
$$T_u=\{P\in T: P_T\le P_u\}\ \ ,$$
and writing 
\begin{equation}\label{upperlowerhalf}
\bigcup_{P\in T}P= (P_T)_u \cup (P_T)_{l} \cup \bigcup_{P\in T_u} P_d \cup \bigcup_{P\in T_d} P_u
\ \ .
\end{equation}
For a convex tree $T$ define the phase plane projection
$$\Pi_{T}=\Pi_{\bigcup P: P\in T}\ \ .$$

We quickly recall some well-known estimates for trees. First, given any tree $T$ and bounded coefficients $\{a_P\}_{P \in T},$ the inequality
\begin{equation} \label{treesingularintegral}
\|\sum_{P \in T_d} a_P \Pi_{P_u} f \|_{p} \leq C_p \sup_{P \in T_d} |a_P| \|f\|_{p}
\end{equation}
holds for $1 < p < \infty$, where $C_p$ depends only on $p$; the proof is by the dyadic Calder\'{o}n-Zygmund method. One then obtains the corresponding square-function estimate
\begin{equation} \label{treesquarefunction}
\|(\sum_{P \in T_d} |\Pi_{P_u} f|^2 )^{1/2}\|_{p} \leq C_p \|f\|_{p}.
\end{equation}
Analogous bounds also hold for $T_u$ (once $\Pi_{P_u}$ is replaced by $\Pi_{P_d}$). Finally, given any convex tree $T$, one can rewrite $\bigcup_{P\in T} P$ as the disjoint union of minimal tiles contained in $\bigcup_{P\in T} P.$ One then immediately sees that $\Pi_Tf$ is pointwise dominated by the dyadic Hardy-Littlewood maximal function of $f$ and hence $\Pi_T$ is bounded from $L^p$ to $L^p$ for $1 < p \leq \infty.$ 

For a tree $T$ and a point $\xi\in \omega_T$ define the enlarged tree $T^{(L)}$ 
to be the set of all bitiles $P$ such that $2^L\xi\in \omega_P$ and
$P$ is contained in a rectangle $2^L P'$ with $P'\in T$. The maximal
element of $T^{(L)}$ is the unique bitile $P$ in $T^{(L)}$ with $I_P=I_T$.
If $T$ is convex then $T^{(L)}$ is also convex. The phase plane projection
$$\Pi_{T^{(L)}}=\Pi_{\bigcup P: P\in T^{(L)}}$$
does not depend on the choice of the frequency $\xi$, because this choice 
is only relevant for the bitiles of $T^{(L)}$  which are contained in $2^LP_T$ 
and these bitiles cover all of $2^LP$ independently of this choice.

Define the trilinear form associated to any subset $\P'\subset \P$  by
$$\Lambda_{\P'}(f_1,f_2,f_3)=
\int \sum_{P\in \P'} \w_{P_d}(x) \Pi_{P_u} f_1(x)
\prod_{j=2}^3 \Pi_{2^LP_d} f_j(x) \, dx\ \ .$$
For a convex collection $\P'$ of tiles define
$$\size(\P',f)=\sup_T |I_T|^{-1/2} \|\Pi_{T} f\|_2$$
or more generally for $L\ge 1$ 
$$\size^{(L)}(\P',f)=\sup_T |I_T|^{-1/2} \|\Pi_{T^{(L)}} f\|_2\ \ ,$$
where in each case the supremum is taken over all convex trees that 
are subset of the collection $\P'$.

\begin{lemma}[{\bf Tree Estimate}]\label{treelemma}
For each $0 < \gamma < 1$ there is a $C_\gamma$ (independent of $L$) such that for any convex tree $T$ and any three bounded functions $f_1$, $f_2$, and $f_3$ in $L^2(\R_+)$ we have
\begin{multline}\label{treeestimate}
|\Lambda_T(f_1,f_2,f_3)|\le C_{\gamma} |I_T| \, \size(T,f_1) \\ \cdot \, \size^{(L)}(T,f_2)^{1-\gamma} \, \size^{(L)}(T,f_3)^{1-\gamma} \|f_2\|^{\gamma}_{{\infty}}\, \|f_3\|_{{\infty}}^{\gamma}.
\end{multline}
\end{lemma}

Proof:
Following the decomposition
(\ref{upperlowertree}) it suffices to prove the estimate for the
three summands of
$$\Lambda_T=\Lambda_{\{P_T\}}+\Lambda_{T_u}+\Lambda_{T_d}$$
separately.
The form $\Lambda_{T_u}$ is estimated by a double application
of Cauchy Schwarz:
$$\Lambda_{T_u}(f_1,f_2,f_3) \le 
\left(\sup_{P\in T} \|\Pi_{P_u} f_1\|_\infty\right)
\prod_{j=2}^3 \left(\sum_{P\in T_u} \|\Pi_{2^LP_d} f_j\|_2^2\right)^{1/2}\ \ .$$
To estimate the first factor on the right-hand-side we consider for any individual bitile 
$P$ the size estimate for the tree $\{P\}$ and obtain
$$\|\Pi_{P_u}f_1\|_\infty \le |I_P|^{-1/2} \|\Pi_{P_u}f_1\|_2\le \size(\{P\},f_1)
\le \size(T,f_1)$$
where the first inequality follows from the fact that $\Pi_{P_u}$ is a rank one
projection onto the space of multiples of $w_{P_u}$.
To estimate the other two factors we observe 
$$\sum_{P\in T_u} \|\Pi_{2^LP_d} f_j\|_2^2\le 
\sum_{P\in (T^{(L)})_u} |\<f_j, w_{P_d}\>|^2
\le \size^{(L)}(T,f_j) |I_T|^{1/2}\ \ .
$$
where the first inequality follows by covering each of the pairwise
disjoint rectangles $2^LP_d$ by tiles of the form $P_d$ with
$P\in T^{(L)}$. Combining these estimates and using the fact that $\size^{(L)}(T,f_j) \leq \|f_j\|_{\infty}$ completes the bound for $T_u$.

The form $\Lambda_{\{P_T\}}$ is estimated similarly, so it 
remains to estimate the form $\Lambda_{T_d}$. 
We have the following identity
\begin{equation}\label{treeprojection}
\Lambda_{T_d}(f_1,f_2,f_3)=\Lambda_{T_d}(\Pi_{T}f_1,\Pi_{T^{(L)}}f_2,\Pi_{T^{(L)}}f_3)
\end{equation}
by an application of the fact shown in \cite{thiele} that for any
two sets $S_1\subset S_2$ in the phase plane which can be written as disjoint unions of tiles 
we have 
\begin{equation}\label{doubleprojection}
\Pi_{S_1}=\Pi_{S_1}\Pi_{S_2}
\end{equation}
It thus remains to show that
\begin{multline*}
|\Lambda_{T_d}(h_1,h_2,h_3)|
\le C|I_T|\, \size(T,f_1) \\ \cdot
\size^{(L)}(T,f_2)^{1-\gamma} \|f_2\|_{\infty}^{\gamma}\, \size^{(L)}(T,f_3)^{1-\gamma}\|f_3\|_{{\infty}}^{\gamma}\ \ ,
\end{multline*}
where $h_1 = \Pi_{T} f_1$ and $h_j = \Pi_{T^{(L)}}f_j$ for $j=2,3$. 

By dilation we may assume $|I_T|=1$. Choose a frequency $\xi\in \omega_T$. 
and define for $l\ge 0$ the interval $\omega_l$ to be the dyadic interval of length
$2^l$ which contains $2^L\xi$. Define
$$\Pi_l:=\Pi_{I_T\times \omega_l}$$ 
and for $l\ge 1$ define
$$\Pi^\Delta_l=\Pi_l-\Pi_{l-1}$$
Then, using a telescoping argument, we may write for $\Lambda_{T_d}(h_1,h_2,h_3)$ 
$$\int \sum_{l=0}^\infty \sum_{P\in T_d: |I_P|=2^{-l}} \w_{P_d}(x) \Pi_{P_u} h_1(x)
\prod_{j=2}^3 \left( \Pi_l h_j(x) + \sum_{m=1}^L \Pi^{\Delta}_{l+m} h_j (x)\right)\, dx$$
The crucial fact then is that for $|I_P|=2^{-l}$ and $m\neq m'$ with at least one of $m$,$m'$ greater
than one we have
\begin{equation}\label{vanishing}
\int \w_{P_d}(x) \Pi_{P_u} h_1(x)
\Pi^{\Delta}_{l+m} h_2(x) \Pi^{\Delta}_{l+m'} h_3(x) \, dx=0\ \ .
\end{equation}
Namely, the product $\w_{P_d}(x) \Pi_{P_u} h_1$ is a multiple of the Haar function on $I_P$.
On the other hand, the product $\Pi^{\Delta}_{l+m} h_2 \Pi^{\Delta}_{l+m'} h_3$ 
has mean zero on either half of $I_P$.
Likewise, (\ref{vanishing}) holds if $m=m'=1$ or if $m=m'=0$ and
$\Pi^{\Delta}_l$ is replaced by $\Pi_l$, because then
the product of the two factors involving $h_2$ and $h_3$ restricted to $I_P$
is the multiple of the square of a Walsh wave packet on this interval and thus constant
on the interval.

Hence we can write $\Lambda_{T_d}(h_1,h_2,h_3)$ as a sum of three terms:
$$\int \sum_{l=0}^\infty \sum_{P\in T_d: |I_P|=2^{-l}} \w_{P_d}(x) \Pi_{P_u} h_1(x)
\Pi_l h_2(x) \Pi^{\Delta}_{l+1} h_3 (x)
\, dx$$
$$+\int \sum_{l= 0}^\infty \sum_{P\in T_d: |I_P|=2^{-l}} \w_{P_d}(x) \Pi_{P_u} h_1(x) 
\Pi^{\Delta}_{l+1} h_2 (x)\Pi_l h_3(x)
\, dx$$
$$+\int \sum_{l= 0}^\infty \sum_{P\in T_d: |I_P|=2^{-l}} \w_{P_d}(x) \Pi_{P_u} h_1(x)
\sum_{m=2}^L  \prod_{j=2}^3  \Pi^{\Delta}_{l+m} h_j (x) \, dx\ \ .$$
The first two summands are estimated as in the case of $T_u$. Namely, the first summand is
\begin{multline*}
\leq  \|\left(\sum_{P\in T_d} |\Pi_{P_u} h_1|^2\right)^{1/2}\|_{{1/\gamma}} \cdot \|\sup_{l\geq 0} |\Pi_{l}h_2|\|_{{2/(1-\gamma)}} 
\\ \cdot \|\left(\sum_{l= 0}^\infty |\Pi_{l+1}^{\Delta} h_3|^2\right)^{1/2}\|_{{2/(1-\gamma)}}.
\end{multline*}
Applying (\ref{treesquarefunction}) to the first factor above, one sees that it is 
\[
\leq C \|\Pi_{T}f_1\|_{{1/\gamma}} \leq C |I_T|^{\gamma} \size(T,f_1)
\]
where the last inequality follows by interpolating $L^2$ and $L^\infty$ estimates for $\gamma < 1/2$ and by applying H\"{o}lder's inequality to the $L^2$ estimate for $\gamma > 1/2$.
Using the $L^p$ boundedness of maximal dyadic averages, the second factor above is
\[
\leq C \|\Pi_{T^{(L)}}f_2\|_{{2/(1-\gamma)}}.
\]
Since $T \subset T$, we have by definition $\|\Pi_{T^{(L)}}f_2\|_{2} \leq |I_T|^{1/2} \size^{(L)}(T,f_2)$, and since $\Pi_{T^{(L)}}$ is bounded on $L^\infty$, we may interpolate to see that the display above is
\[
\leq C |I_T|^{(1-\gamma)/2} \size^{(L)}(T,f_2)^{1-\gamma} \|f_2\|_{\infty}^{\gamma}.
\]
Similarly, but using the square function bound again, the third factor is 
\[
\leq C |I_T|^{(1-\gamma)/2} \size^{(L)}(T,f_3)^{1-\gamma} \|f_3\|_{\infty}^{\gamma}, 
\]
and we thus obtain the desired estimate for the product.

To bound the third summand, we change the order of summation and then estimate
by H\"{o}lder:
$$ |\int \sum_{l=2}^\infty 
\left(\sum_{P\in T_d: 2^{-(l-2)}\le |I_P|\le 2^{-(l-L)}} \w_{P_d}(x) \Pi_{P_u} h_1(x)\right)
\prod_{j=2}^3  \Pi^{\Delta}_{l} h_j (x) \, dx|$$
\begin{equation}\label{twofactors}
\le \left\|\sup_l |\sum_{P\in T_d: 2^{-(l-2)}\le |I_P|\le 2^{-(l-L)}} \w_{P_d} \Pi_{P_u} h_1| \right\|_{{1/\gamma}}
\left\|\sum_{l=2}^{\infty} 
\prod_{j=2}^3  
|\Pi^{\Delta}_{l} h_j | \right\|_{{1/(1-\gamma)}}
\end{equation}
For the second factor we have
$$\left\|\sum_{l=2}^{\infty} 
\prod_{j=2}^3  
|\Pi^{\Delta}_{l} h_j | \right\|_{{1/(1-\gamma)}}
\le \prod_{j=2}^3 \|\left(\sum_{l=2}^{\infty}  |\Pi_{l}^{\Delta} h_j|^2\right)^{1/2}\|_{{2/(1-\gamma)}} 
\ .$$
which, as before, is 
\[
\leq C |I_T|^{(1-\gamma)} \size^{(L)}(T,f_2)^{1-\gamma} \|f_2\|_{\infty}^{\gamma} \size^{(L)}(T,f_3)^{1-\gamma} \|f_3\|_{\infty}^{\gamma}\ \ .
\]

To estimate the first factor in (\ref{twofactors}) we observe that
by the triangle inequality it suffices to estimate
$$\left\|\sup_l |\sum_{P\in T_d: |I_P|> 2^{l}} \w_{P_d} \Pi_{P_u} h_1| \right\|_{{1/\gamma}}\ \ .$$
Since $\w_{P_d} \Pi_{P_u} f_1$ is a linear combination of Haar functions at level $|I_P|$,
the truncation operator to $|I_P|> 2^{l}$ can be replaced by an averaging operator to dyadic 
intervals of length $2^{l}$. Thus, by the Hardy Littlewood maximal theorem, the display above is
\[
\left\|\sum_{P\in T_d} \w_{P_d} \Pi_{P_u} h_1 \right\|_{{1/\gamma}}\ \ .
\]
Since $|\w_{P_d}| \leq 1$ for each $P$, we obtain from (\ref{treesingularintegral}) that the display above is
$$\leq C \left\|h_1\right\|_{{1/\gamma}}\le C |I_T|^{\gamma}\size(T, f_1)\ \ .$$
Combining the estimates for the two factors of (\ref{twofactors}) proves the desired
bound for complete trees and ends the proof of Lemma \ref{treelemma}.

\begin{lemma}[{\bf Tree Selection}]\label{treeselection}
Assume $\P$ is a convex collection of bitiles contained
in the strip $\R_+\times [0,2^N)$
with $\size^{(L)}(\P,f)\le 2^{-k}$. Then we can write
$\P$ as the union of a convex set of bitiles $\P'$ and a collection $\T$
of convex trees such that
\begin{equation}\label{l1N}
\|\sum_{T\in \T} 1_{I_T}\|_1 \le C 2^{2k}\|f\|_2^2\ \ ,
\end{equation}
\begin{equation}\label{bmoN}
\|\sum_{T\in \T} 1_{I_T}\|_{\rm BMO} \le C 2^{2k}\|f\|_\infty^2\ \ ,
\end{equation}
with constants independent of $N$ and $L$ and 
$\size^{(L)}(\P',f)\le 2^{-k-1}$.
\end{lemma}
Note that the case $L=0$ corresponds to a statement for $\size(\P,f)$.

Proof: By scaling it suffices to prove the lemma for $k=0$.
We shall first see that we may reduce to the case that all bitiles $P\in \P$
satisfy
\begin{equation}\label{lsizetile}
\size^{(L)}(\{P\},f)\le 2^{-4}\ \ .
\end{equation}
If there is a bitile $P\in \P$ which violates (\ref{lsizetile}), 
then we pick one such bitile $P_1$ which maximizes
$I_{P_1}$ and set $T_1=\{P'\in \P: P'\le P_1\}$ and $\P_1=\P\setminus T_1$.
Both $T_1$ and $\P_1$ are convex.
Then we iterate this procedure with $\P_1$, provided there
is a bitile in $\P_1$ which violates (\ref{lsizetile}), and so on.
The selected bitiles are all pairwise disjoint. For assume not, 
then $P_m\le P_k$ for some $m,k$
and by choice of these bitiles we necessarily have $k<m$. But then
$P_m$ should have been in the tree $T_k$ and would not have
been available for selection at the $m$-th step, a contradiction.
By vertical dilation, the rectangles $2^LP_m$ with $P_m$ a selected bitile
are pairwise disjoint. We therefore have
$$\sum_{k=1}^n \|\Pi_{2^LP_k} f\|_2^2\le  \|f\|_2^2$$
and hence
$$\sum_{k=1}^n |I_{T_k}|\le 2^{8}\|f\|_2^2$$
which proves (\ref{l1N}) for the set of selected trees. 
Moreover, for every
dyadic interval $I$ we have
$$\sum_{k: I_{P_k}\subset I} \|\Pi_{2^LP_k} f\|_2^2\le \|f1_I\|_2^2\ \ .$$
This gives (\ref{bmoN}) for the set of selected trees.
Since all bitiles in the collection $\P$ satisfy $|I_P|\ge 2^{-N}$,
estimate (\ref{l1N}) shows that the selection process must have stopped
after finitely many steps, and the remaining collection
has no bitiles violating (\ref{lsizetile}).
For the rest of the argument we assume that
all bitiles in $\P$ satisfy (\ref{lsizetile}).

By (\ref{upperlowerhalf}) it suffices to show 
that for the collection $\P'$ we are about to construct we have
for every tree $T \subset \P'$ 
\begin{equation}\label{upperest}
\sum_{P\in T_u} \|\Pi_{2^LP_d} f\|_2^2\le 2^{-4}|I_T|
\end{equation}
\begin{equation}\label{lowerest}
\sum_{P\in T_d} \|\Pi_{2^LP_u}f\|_2^2\le 2^{-4}|I_T|
\end{equation}
because we already have 
$$\|\Pi_{2^LP_T} f\|_2^2\le 2^{-8} |I_T|\ \ .$$
We will prove that
one can take away a collection $\T$ of trees satifying (\ref{l1N}) and (\ref{bmoN})
such that for the remaining collection $\P'$ of bitiles we have 
the bound (\ref{upperest}) for all convex trees $T\subset \P'$. Inequality (\ref{lowerest}) is covered by an analogous argument which will be omitted. 

To do so, we again iteratively select trees. If there is a
tree in the collection $\P$ which violates (\ref{upperest}), we choose
one such tree $S_1$ with maximal element $P_1$ such that the left endpoint 
of $\omega_{P_1}$ is minimal. We may assume $f$ is non-zero and has finite $L^2$ norm,
hence violation of (\ref{upperest}) implies an upper bound for $I_{S_1}$ and thus there
are only finitely many possible choices for the interval $\omega_{P_1}$
contained in $[0,2^N)$ and thus one of the choices attains the minimum for the
left endpoint of $\omega_{P_1}$.

Then we define $T_1=\{P\in\P: P\le P_1\}$ and 
$\P_1=\P\setminus T_1$. Both $T_1$ and $\P_1$ are convex and $T_1$ contains $S_1$.
Then we iterate this procedure as long as the remaining collection $\P_n$ contains 
a tree which violates (\ref{upperest}).
We prove (\ref{l1N}) and (\ref{bmoN}) for the collection
of selected trees.

Let ${S}_k$ and ${S}_m$ be two different selected trees. We claim that if
$P_k$ and $P_m$ are bitiles in the respective trees, then we have
that $(P_k)_d$ and $(P_m)_d$ are disjoint. For assume not, then
without loss of generality 
$$(P_k)_d\le (P_m)_d\  .$$
Then also $(P_k)_d\le (P_m)_u$ because $P_k$ and $P_m$ have to be different.
This implies that $m<k$ by the choice of trees. 
Then $P_k$ should have been selected for the tree $T_m$ and should not
have been available for $S_k$. This is the desired contradiction and establishes
that $(P_k)_d$ and $(P_m)_d$ are disjoint. Hence we have for the collection of selected trees
$$\sum_{m}\left(\sum_{P\in \tilde{S}_m} \|\Pi_{2^LP_d}f\|_2^2\right) \le \|f\|_2^2\ \ .$$
Combining this with the violation of (\ref{upperest}) shows the desired estimate (\ref{l1N}). 
Estimate (\ref{bmoN}) then follows again by localization.

By virtue of (\ref{l1N}) and the lower bound $2^{-N}\le |I_T|$ we are guaranteed that 
the tree selection process stops and the remaining collection does not contain a tree
that violates (\ref{upperest}).
This completes the proof of Lemma \ref{treeselection}\ \ .


\section{Proof of the Main Theorems} \label{proof}

We first prove Theorem \ref{main1} in the open triangle $c$
of Figure \ref{hexagon}. Then we prove certain restricted weak type 
bounds in the open diamond $b_3\cup d_{23}$
and use interpolation techniques to obtain Theorem \ref{main1}
in the open triangle $b_3$ and Theorem \ref{main2} in the
open triangle $d_{23}$. Symmetric arguments can be applied to the 
diamonds $b_1\cup d_{21}$, $b_1\cup d_{31}$ and $b_2\cup d_{32}$.
The argument is however not entirely symmetric, it does not apply to
the forbidden diamonds $b_2\cup d_{1,2}$ or $b_3\cup d_{13}$.
The full extent of Theorems \ref{main1} and \ref{main2} then
follows by interpolation.

\begin{proposition}[{\bf Triangle $c$}]\label{regionc}
For each $i=1,2,3$ let  $E_i$ be a subset of $\R_+$ and 
$f_i$ a measurable function bounded by the characteristic function
of $E_i$. Let $0 \le \alpha_1 \leq 1/2$, $0 \le \alpha_i < 1/2$ for $i = 2,3$, and assume
$\alpha_1+\alpha_2+\alpha_3=1$.
Let $j$ be an index such that $|E_j|$ is maximal. Then there is a 
major subset ${E_j}'$ of $E_j$ depending
on $E_1,E_2,E_3$ such that if $f_j$ is also supported in ${E_j}'$
we have
$$\Lambda_L(f_1,f_2,f_3)\le C_{\alpha_1,\alpha_2,\alpha_3}
|E_1|^{\alpha_1}|E_2|^{\alpha_2}|E_3|^{\alpha_3}$$
with a constant $C_{\alpha_1,\alpha_2,\alpha_3}$ independent of $L$. If the three sets $E_i$,
have measure within a factor of four of each other, we may 
choose ${E_j}'=E_j$.
\end{proposition}

Proof:
Dilating by a power of $2$ we may assume $1\le |E_j|<2$. Define the exceptional set
$$F=\bigcup_{i\neq j} \{x: M_2(1_{E_i}/|E_i|^{1/2})(x)>2^{10}\}\ \ .$$
By the Hardy Littlewood maximal theorem the measure of $F$
is less than one half and we may define the major subset ${E_j}'=E_j\setminus F$.
The set $F$ is empty if the measure of all $E_i$ is at least one fourth.

Given functions $f_i$ as in the proposition define the normalized functions
$g_i=f_i |E_i|^{-1/2}$ for $i=1,2,3$. Fixing $\gamma$ with 
$$0 < \gamma < \min_{\alpha = 1/4,\alpha_2, \alpha_3} 1 - 2\alpha\ \ ,$$ 
it suffices to show
$$\Lambda_L(g_1,g_2,g_3)\le |E_2|^{-\gamma/2} |E_3|^{-\gamma/2}\ \ .$$
Let $\P$ be the convex set of all bitiles $P$ in the strip $\R_+\times [0,2^N)$ 
such that $I_P$ is not contained in $F$. Since $g_j$ vanishes on $F$ we have 
$$\Lambda_{\P}(g_1,g_2,g_3)=\Lambda_L(g_1,g_2,g_3)\ \ .$$
Outside the set $F$, the $M_2$ maximal function of each of the three
functions $g_i$ is bounded by a universal constant, hence we have
$$\size(\P,g_i)\le C,\ \ \size^{(L)}(\P,g_i)\le C\ \ \ .
$$
Applying Lemma \ref{treeselection} repeatedly, we define a decreasing nested sequence
of convex subsets ${\P}_k$ of the set ${\P}_0:=\P$ such that
$$\size({\P}_k,g_i)\le C2^{-k},\ \ \size^{(L)}({\P}_k,g_i)\le C2^{-k}$$
for each $i$ and ${\P}_{k-1}$ is the disjoint union of ${\P}_{k}$ and a collection 
${\T}_k$ of convex trees with 
$$\sum_{T\in {\T}_k} |I_T| \le C 2^{2k}\ \ .$$
We then have
\begin{align*}
\Lambda_\P(g_1,g_2,g_3) &=\sum_{k\ge 1} \sum_{T\in {\T}_k} \Lambda_T(g_1,g_2,g_3) \\
&\le \sum_{k\ge 1} \sum_{T\in {\T}_k} C|I_T| \size(T,g_1)\prod_{i = 2,3} \size^{(L)}(T,g_i)^{1-\gamma} |E_i|^{-\gamma/2} \\
&\le \sum_{k\ge 1} \sum_{T\in {\T}_k} C 2^{-3k +2\gamma} |I_T| |E_2|^{-\gamma/2}|E_3|^{-\gamma/2} \\
&\le \sum_{k\ge 1} C 2^{-k(1 - 2\gamma)} |E_2|^{-\gamma/2}|E_3|^{-\gamma/2} \\
&\le C |E_2|^{-\gamma/2}|E_3|^{-\gamma/2}.
\end{align*} 
This completes the proof of the proposition.

By interpolation as in \cite{mtt0} this proposition proves Theorem \ref{main1}
in the region $2<p_j<\infty$. First one proves that the proposition holds 
in this region for $f_j$ not necessarily supported on ${E_j}'$. Namely one splits 
$f_j$ into $f_j1_{{E_j}'}+f_j1_{E_j\setminus {E_j}'}$. On the first summand the conclusion
of the proposition gives the desired bound, while on the second summand one 
iterates the proposition with $E_j$ replaced by $E_j\setminus {E_j}'$.
One continues the iteration until all three sets are of comparable size at 
which time the proposition holds already with the major subset being the full set. 
The various estimates throughout the iteration process are summable provided 
$\alpha_j < 1/2$ for all $j$. With this variant of the proposition established, 
one applies standard multilinear Marcinkiewicz interpolation to obtain the strong type estimate 
in the region $2<p_j<\infty$.

\begin{proposition}[{\bf Diamond $b_3\cup d_{23}$}]\label{diamond}
Let $0<\epsilon,\alpha<1/2$.
For each $i=1,2,3$ let $E_i$ be a measurable subset of $\R_+$ and 
let $f_i$ be a measurable function bounded by the characteristic function
of $E_i$.
Assume
$|E_3|<|E_2|$ and $1\le |E_2|\le 2$.
Then there is a major subset 
${E_2}'$ of $E_2$ depending on $E_1,E_2,E_3$ 
such that if $f_2$ is supported in ${E_2}'$ we have
$$\Lambda_L(f_1,f_2,f_3)\le C_{\alpha,\epsilon}
|E_1|^{\alpha}
|E_3|^{1-\epsilon}
$$
with a constant $C_{\alpha,\epsilon}$ independent of $L$. 
\end{proposition}

Proof: 
Define the exceptional set $F$ to be
$$
\{x: M_2(1_{E_1}/|E_1|^{\alpha})(x)>2^{10}\}
\cup
\{ M_{1/(1-\epsilon)} (1_{E_3}/|E_3|^{1-\epsilon})(x)>2^{10}\}
\ \ .$$
We may assume that $|E_3|$ is sufficiently small so that $E_3$
is contained in $F$, or else the desired estimate is trivial from
the already established case of Theorem \ref{main1} 
in the vicinity of $\alpha_2=1/2$ and $\alpha_3=1/2-\alpha$.

By the Hardy Littlewood maximal theorem the measure of $F$
is less than $1/2$ and we may define the major subset ${E_2}'=E_2\setminus F$.
Given a triple of functions as in the proposition define the normalized
functions 
$$g_1=f_1 |E_1|^{-\alpha},\ \ g_2=f_2,\ \ g_3=f_3 |E_3|^{\epsilon-1}\ \ .$$
Let $\P$ be the convex collection of all bitiles in the strip
$\R_+\times [0,2^N)$ such that $I_P$ is not contained in $F$. We have
$$\size(\P,g_1)\le C \ \ .$$
Applying Lemma \ref{treeselection} repeatedly 
we obtain a nested sequence 
${\P}_k$ of convex subsets of $\P={\P}_0$ 
such that we have
$$\size(\P_k,g_1)\le C2^{-k}$$
and, interpolating between (\ref{bmoN}) and (\ref{l1N}),
the set $\P_{k-1}$ is the disjoint union of $\P_k$ and a collection
$\T_k$ of convex trees with
\begin{equation}\label{pcounting}
\|\sum_{T\in {\T}_k} 1_{I_T}\|_p \le C 2^{2k}
\end{equation}
for $1/p=2\alpha$.
We shall fix a $k\ge 1$ and prove
\begin{equation}\label{fixedk}
|\sum_{T\in \T_k} \Lambda_T(g_1,g_2,g_3)|\le C 2^{-\epsilon k} 
\ \ ,
\end{equation}
which will clearly finish the proof of the proposition.

For $P\in \P_k$ define $\p_P$ to be the
set of all minimal tiles (those with spatial interval of length $2^{-L}|I_P|$) contained in 
$2^LP_d$ for which $I_p$ is not contained in $F$. Then we may write
$$ \int \w_{P_d}(x) \Pi_{P_u} g_1(x)
\prod_{j=2}^3 \Pi_{2^LP_d} g_j(x) \, dx$$
$$ = \int \w_{P_d}(x) \Pi_{P_u} g_1(x)
\sum_{p\in \p_P} \prod_{j=2}^3 
\< g_j,w_p\> w_p(x) \, dx\ \ .$$
because the wave packets $w_p$ are disjointly supported as $p$ runs through $\p_P$
and $\<g_2,w_p\>=0$ if $I_p\subset F$. Note that in this argument the symmetry between the three 
indices $1,2,3$ is broken. The role played by the index $2$ in this argument, namely the use of
vanishing of $g_2$ on $F$, could symmetrically be taken by the index $3$, but not by the 
index $1$. This is the only place in the proof of this proposition, where the symmetry is 
broken in an essential way.

Let $\I$ be the collection of maximal dyadic intervals contained in the set 
$F$. For an interval $I\in \I$ let $\p_I$ be the collection of tiles $p$
with time interval $I$ which intersect a tile $p'$ in $\p_P$ for some
$P\in \bigcup_{T\in \T_k} T$. We have $p\le p'$ for such tiles, and the relation is strict
in the sense $p'\neq p$. Then we have also $p \leq 2^LP_T$ and hence there is at most one
element in $\p_I$ which intersects with a given tree in $\T_k$.
Hence $\p_I$ has at most $N_I$ elements where $N_I$ is the constant value of
the function 
$$\sum_{T\in \T_k} 1_{I_T}$$
on the interval $I$.

Let $a_{I}$ be the orthogonal projection of $g_3$ onto the span
of wave packets associated to the tiles in $\p_I$, and let $a=\sum a_I$.
Since $g_3$ is supported on the union of the intervals $I\in \I$, we have for every 
tree $T\in \T_k$
$$\Lambda_{T}(g_1,g_2,g_3)=\Lambda_T(g_1,g_2,a)\ \ .$$
We have (H\"older and Hausdorff Young on $I$), 
$$\|a_{I}\|_2^2= \sum_{p\in \p_I}|\<g_3,w_p\>|^2$$
$$\le N_I^{1-2\epsilon} (\sum_{p\in \p_I}|\<g_3,w_p\>|^{1/\epsilon})^{2\epsilon}$$
$$\le N_I^{1-2\epsilon} \|1_Ig_3\|_{1/(1-\epsilon)}^2 |I|^{2\epsilon-1}\le C N_I^{1-2\epsilon}|I|$$
where in the  second to last inequality we interpolate between $L^2 \rightarrow \ell^2$ and $L^1 \rightarrow \ell^{\infty}$ bounds and in the last inequality we used the bound for the $M_{1/(1-\epsilon)}$ - function of
$g_3$ at some point of the dyadic parent of $I$.
Hence we have
$$\|a\|_2^2\le  C \sum_I N_I^{1-2\epsilon} |I|$$
$$\le C (\sum_I N_I^p |I|)^{(1-2\epsilon)/p}(\sum_I |I|)^{1-(1-2\epsilon)/p}\le C 2^{2k(1-2\epsilon)}\ \ .$$
Also note that $\|a_I\|_{\infty} \leq C N_I.$

Define $\I_{0}$ to be the set of all intervals in
$\I$ such that $N_I$ is at most $2^{2k}$. For $m>0$ define
$\I_{m}$ to be the set of all intervals in $\I$
such that $N_I$ is between $2^{2(k+m-1)}$ and $2^{2(k+m)}$.
Split
$$a=\sum_{m\ge 0} a_m$$
accordingly, i.e., $a_m$ is supported on the union of intervals
in $\I_{m}$. We have 
$$\|a_0\|_2^2\le C 2^{2k(1-2\epsilon)}$$
and for $m>0$
$$\|a_m\|_2^2\le C \sum_{I\in \I_{m}} 2^{2(k+m)(1-2\epsilon)} |I|$$
$$\le  C2^{2(k+m)(1-2\epsilon)} 2^{- 2p(k+m)} \|\sum_I N_I 1_I\|_p^p$$
$$\le C 2^{2(k+m)(1-2\epsilon)} 2^{-2pm}\ \ .$$
Moreover, we have for every tree $T$ contained in $\bigcup_{T\in \T_k}T$
$$\|a_m\|_{L^2(I_T)}^2\le C \sum_{I\in \I_{m}:I\subset I_T} 2^{2(k+m)(1-2\epsilon)} |I|
\le C 2^{2(k+m)(1-2\epsilon)} |I_T|
$$
and hence 
$$\size^{(L)}(\bigcup_{T\in \T_k} T, a_m)\le C 2^{(k+m)(1-2\epsilon)}\ \ .$$
We shall fix $m$ and prove
\begin{equation}\label{fixedkm}
|\sum_{T\in \T_k}\Lambda_k(g_1,g_2,a_m)|\le C  2^{-\epsilon k} 2^{- \epsilon m}\ , 
\end{equation}
which will imply (\ref{fixedk}) and finish the proof of the proposition.
Normalize 
$$\tilde{a}_m= a_m 2^{-(k+m)(1-2\epsilon)} 2^{pm}$$ 
so that $\|\tilde{a}_m\|_2\le C$, then it clearly suffices to prove
$$|\sum_{T\in \T_k}\Lambda_k(g_1,g_2,\tilde{a}_m)|\le  C  2^{- k(1-\epsilon)} 2^{m\epsilon}\ . $$

With Lemma \ref{treeselection} we decompose
$\bigcup_{T\in \T_k}T$ into collections $\bigcup_{T\in \T_{k,l}} T$ with $-pm\le l< pk$
and a remainder set $\P_{k,pk}$ 
such that each tree in $\T_{k,l}$ satisfies
$$\size^{(L)}(T,g_2)\le C 2^{-l},\ \ \size^{(L)}(T,\tilde{a}_m)\le C 2^{-l}$$
and we have
$$\sum_{T\in \T_{k,l}} |I_T|\le 2^{2l} $$
and we have for every tree in the remainder set $\P_{k,pk}$
$$\size^{(L)}(T,g_2)\le C 2^{-pk},\ \ \size^{(L)}(T,\tilde{a}_m)\le C 2^{-pk}\ \ .$$
Note that we also have 
$$\size^{(L)}(T,g_2)\le C$$
for every selected tree since $g_2$ is bounded by a universal constant.

Applying Lemma \ref{treelemma} with a small value of $\gamma$ to be determined later, we obtain
$$|\sum_{T\in \T_k} \Lambda_{T_k} (g_1,g_2,\tilde{a}_m)|$$
$$\le C \sum_{-pm \le l <pk} 
2^{-k} \min(1,2^{-l})^{1-\gamma} 2^{-l(1-\gamma)} 2^{2l} 2^{\gamma(k+m)} $$
$$+|\Lambda_{\P_{k,pk}}(g_1,g_2,\tilde{a}_m)|$$ 
The first term is a converging geometric series for $l<0$ and a diverging geometric series
for $l>0$ and hence 
$$C \sum_{-pm \le l <pk}  2^{-k} \min(1,2^{-l})^{(1-\gamma)} 2^{-l(1-\gamma)} 2^{2l} 2^{\gamma(k+m)} \le C 2^{- k} 2^{2p\gamma k} 2^{\gamma(k+m)}.$$
To estimate the term involving $\P_{k,pk}$ we use that $\P_{k,pk}$ is a subset of the union
of trees in $\T_k$. For every tree $T'$ in $\T_k$ we can write $\P_{k,pk}\cap T'$
as a union of trees $T$ with $\sum |I_T|\le |I_{T'}|$. By (\ref{pcounting}) we can write
$\P_{k,pk}$ as the union of trees with
$$\sum_{T}|I_T|\le 2^{2pk}\ \ .$$
Hence we have
$$|\Lambda_{\P_{k,pk}}(g_1,g_2,\tilde{a}_m)|\le C2^{-k} 2^{-pk(1-\gamma)} 2^{-pk(1-\gamma)} 2^{2pk} 2^{\gamma(k+m)}.$$
Choosing $\gamma < \epsilon/(1 + 2p)$, both terms are controlled by  $2^{- k(1-\epsilon)} 2^{m\epsilon}$ as desired.
This proves (\ref{fixedkm}) and completes the proof of the proposition.

By scaling we can remove the restriction $1\le |E_2|\le 2$.
This proves Theorem \ref{main2} in the open triangle $d_{23}$.
In the open triangle $b_3$, we may iterate Proposition \ref{diamond}
as in the discussion of the previous proposition, and as a result obtain that 
the proposition holds for ${E_2}'=E_2$. Then Theorem \ref{main1} holds in the 
open triangle $b_3$ by multilinear Marcinkiewicz interpolation.


\end{document}